\documentclass[12pt]{article}
\usepackage{amssymb, amsmath}
\begin{document}
\renewcommand{\theequation}{\arabic{section}.\arabic{equation}}
\newtheorem{theorem}{Theorem}[section]
\newtheorem{lemma}{Lemma}[section]
\newtheorem{definition}{Definition}
\newtheorem{pro}{Proposition}[section]
\newtheorem{cor}{Corollary}[section]
\newcommand{\n}{\nonumber}
\newcommand{\tv}{\tilde{v}}
\newcommand{\tw}{\tilde{\omega}}
\renewcommand{\t}{\theta}
\newcommand{\w}{\omega}
\newcommand{\e}{\varepsilon}
\renewcommand{\a}{\alpha}
\renewcommand{\l}{\lambda}
\newcommand{\vare}{\varepsilon}
\newcommand{\s}{\sigma}
\renewcommand{\o}{\omega}
\renewcommand{\O}{\Omega}
\newcommand{\bb}{\begin{equation}}
\newcommand{\ee}{\end{equation}}
\newcommand{\bq}{\begin{eqnarray}}
\newcommand{\eq}{\end{eqnarray}}
\newcommand{\bqn}{\begin{eqnarray*}}
\newcommand{\eqn}{\end{eqnarray*}}
\title{Nonexistence of  self-similar singularities in
the ideal magnetohydrodynamics}
\author{Dongho Chae\thanks{This work was supported partially by  KRF Grant(MOEHRD, Basic
Research Promotion Fund) and the KOSEF Grant no.
R01-2005-000-10077-0. Part of the work was done, while the author
was visiting RIMS, Kyoto University. He would like to  thank
Professor Hisashi Okamoto for his hospitality during the visit and
useful discussions.\newline
 }\\
Department of Mathematics\\
              Sungkyunkwan University\\
               Suwon 440-746, Korea\\
              {\it e-mail : chae@skku.edu}}
 \date{}
\maketitle
\begin{abstract}
In this paper we exclude the scenario of apparition of finite time
singularity in the form of self-similar singularities in the ideal
magnetohydrodynamic equations, assuming  suitable integrability
conditions on the vorticity and the magnetic field. We also consider
more sophisticated possibility of asymptotically self-similar
singularities, which means that the local classical solution
converges to the self-similar profile as we approaches to the
possible time of singularity. The scenario of asymptotically
self-similar singularity is also excluded under suitable conditions
on the profile. In the 2D magnetohydrodynamics the magnetic field
evolution equations reduce to a divergence free transport equation
for a scalar stream function. This helps us to improve the above
nonexistence theorems on the self-similar singularities, in the
sense that we only need weaker integrability conditions on the
profile to prove the results.
\end{abstract}
\section*{Introduction}
We are  concerned  on the possibility of finite time singularity in
the Cauchy problem of the ideal magnetohydrodynamic equations in
$\Bbb R^n$, $n=2,3$.

\[
\mathrm{ (MHD)}
 \left\{ \aligned
 &\frac{\partial v}{\partial t} +(v\cdot \nabla )v =(b\cdot\nabla)b-\nabla (p +\frac12 |b|^2), \\
 &\frac{\partial b}{\partial t} +(v\cdot \nabla )b =(b \cdot \nabla )
 v,\\
 &\quad \textrm{div }\, v =\textrm{div }\, b= 0 ,\\
  &v(x,0)=v_0 (x), \quad b(x,0)=b_0 (x)
  \endaligned
  \right.
  \]
where $v=(v_1, \cdots , v_n )$, $v_j =v_j (x, t)$, $j=1,\cdots,n$,
is the velocity of the flow, $p=p(x,t)$ is the scalar pressure,
$b=(b_1, \cdots , b_n )$, $b_j =b_j (x, t)$, is the magnetic field,
and $v_0$, $b_0$ are the given initial velocity and magnetic field,
 satisfying div $v_0 =\mathrm{div}\, b_0= 0$, respectively.
 The system (MHD) is the incompressible Euler equations coupled with the
 magnetic field evolution equations, and the question of finite time
 singularity/global regularity is an outstanding open problem
 in the mathematical fluid mechanics, similarly to the case of the 3D
 Euler equations(see \cite{cha3}, and references therein).
 There are some numerical approaches to this problem(see e.g. \cite{gra, gib} and
 references therein).  The blow-up criterion similar to the
 Beale-Kato-Majda's one for incompressible Euler system(\cite{bea})
  is obtained by Caflisch, Klapper and Steele(\cite{caf},
 see also \cite{ohk}). In this paper our aim is to  consider the
 possibility of finite time apparition of singularity for (MHD) in the form of
  self-similar
 type, and exclude it. The organization of this paper is the
 following:\\
 In the section 1 we first exclude  self-similar singularities of
 the Leray type for the Navier-Stokes
 equations introduced in \cite{ler}. In the section 2, using the result of section 1, we
 exclude more sophisticated self-similar singularity, called the asymptotically self-similar
 singularities. In the section 3  we present improved theorems of
 sections 1 and 2 in the case of planar magnetohydrodynamics.

\section{Self-similar singularities}
 \setcounter{equation}{0}

In order to be more precise on the notion of self-similar
 singularities we begin by the following observation on the scaling property of (MHD):
  if $(v, b, p)$ is a
solution of  (MHD) corresponding to the initial data $(v_0, b_0)$,
then for any $\lambda
>0$ and $\alpha \in \Bbb R $ the functions
$$
  v^{\lambda, \alpha}(x,t)=\lambda ^\alpha v (\lambda x, \l^{\a +1}
  t), \quad  b^{\lambda, \alpha}(x,t)=\lambda ^\alpha b (\lambda x, \l^{\a +1}
  t)$$
 and
  $$ \quad p^{\l, \a}(x,t)=\l^{2\a}p(\l x, \l^{\a+1} t )
$$
  are also solutions with the initial data
  $ v^{\lambda, \alpha}_0(x)=\lambda ^\alpha v_0
   (\lambda x)$,  $ b^{\lambda, \alpha}_0(x)=\lambda ^\alpha b_0
   (\lambda x)$.
 In view of the above scaling
  property the  self-similar blowing up
  solution $(v(x,t), b(x,t))$ of the system (MHD) should be of the form,
  \bq
  \label{vel}
 v(x, t)&=&\frac{1}{(T_*-t)^{\frac{\a}{\a+1}}}
V\left(\frac{x}{(T_*-t)^{\frac{1}{\a+1}}}\right),\\
  \label{mag}
 b(x, t)&=&\frac{1}{(T_*-t)^{\frac{\a}{\a+1}}}
B\left(\frac{x}{(T_*-t)^{\frac{1}{\a+1}}}\right),\\
\label{pre}
  p(x,t)&=&\frac{1}{(T_*-t)^{\frac{2\a}{\a+1}}}
P\left(\frac{x}{(T_*-t)^{\frac{1}{\a+1}}}\right)
 \eq
 for  $\a \neq -1$ and $t$
 close to the possible blow-up time $T_*$. If we substitute (\ref{vel})-(\ref{pre})
 into (MHD), then we find that $(V,B,P)$ should be a solution of the
 stationary system:
\bb\label{mhdleray}
 \left\{
 \aligned
  &\frac{\a}{\a+1} V+\frac{1}{\a +1} (y\cdot \nabla)
V+(V\cdot \nabla)
 V=(B\cdot \nabla )B -\nabla (P +\frac12 |B|^2)\\
&\frac{\a}{\a+1} B+\frac{1}{\a +1} (y\cdot \nabla)B+(V\cdot
 \nabla )B=(B\cdot \nabla )V\\
 &\qquad\mathrm{ div}\, V=\mathrm{div}\, B=0
 \endaligned
 \right.
 \ee
 Conversely, if $(V, B, P)$ is a smooth solution of the system (\ref{mhdleray}),
 then the triple of functions $(v,b,p)$ defined  by (\ref{vel})-(\ref{pre}) is
 a smooth solution of (MHD) for $t\in (0, T_*)$, which blows up at $t=T_*$.
 The search for self-similar singularities of the form  similar to (\ref{vel})-(\ref{pre})
 (more precisely with $\a =1$ and $B=b=0$) was suggested by Leray  for
 the 3D Navier-Stokes
 equations in \cite{ler}, and the possibility was first excluded by Ne$\check{c}$as,
 Ru$\check{z}$i$\check{c}$ka
  and $\check{S}$ver$\acute{a}$k in \cite{nec} under the condition of $V\in L^3(\Bbb R^3)\cap H^1_{loc}
 (\Bbb R^3 )$, the result of which was generalized later by Tsai in
 \cite{tsa}. Their proofs crucially depends on the maximum principle
 of the Leray system,
$$
\frac12 V+\frac12 (y\cdot \nabla) V+(V\cdot \nabla)
 V=-\nabla P +\Delta V, \quad \mathrm{div}\,V=0,
 $$
which corresponds to the Navier-Stokes version of the
 system (\ref{mhdleray}). The maximum principle, in turn, is possible due to the
 dissipation term, $\Delta V$,
 in the Leray system, which is originated from the dissipation term of the
 Navier-Stokes equations. Due to this fact there was difficulty in extending the
 nonexistence results for the self-similar singularity
 of the 3D Navier-Stokes system to the 3D
 Euler equations, applying  the similar method to \cite{nec} or \cite{tsa}. Recently,
 the author of this paper discovered completely
 new argument to prove nonexistence of the self-similar singularity of
 the 3D Euler system under suitable integrability condition on the
 vorticity(\cite{cha1}). In this section we apply that method iteratively
  to prove nonexistence of self-similar singularity for (MHD).
 The precise theorems are
 stated and proved in the subsections below.

\subsection{Statement of the theorems}

 Below we denote by $C_0 ^m (\Bbb R^n)$ the collection of
$C^m (\Bbb R^n)$ functions with all the derivatives up to the $m-$th
order vanishing at infinity. The class of $C^m (\Bbb R^n )$
functions with compact support will be denoted by $C^m_c (\Bbb
R^n)$.
\begin{theorem}
Suppose there exists $T_*>0$ and  $\a \neq -1$ such that we have  a
 representation of a solution $(v, b)$
  to (MHD) by
  (\ref{vel})-(\ref{mag})  for all $t \in (0 ,T_*)$   with $(V, B)$
  satisfying the following
  conditions:
\begin{itemize}
 \item[(i)] $(V, B)\in [C_0 ^1 (\Bbb R^n)]^2$.
  \item[(ii)] There exists $q_1 >0$ such that
 $(\Omega,  B)\in [L^q (\Bbb R^3
  )]^2$ for all $q\in (0, q_1)$, where  $\O=\mathrm{curl}\, V$.
\end{itemize}
  Then, $V=B=0$.
\end{theorem}\noindent{\textsf{Remark }1.1} Due to the  condition (i) on
$V$ we  can exclude the possibility of the self-similar singular
solution $(v,b,p)$ of the form,
 \bq\label{example}
 v(x, t)&=&\frac{1}{(T_*-t)^{\frac{\a}{\a+1}}}
\nabla h \left(\frac{x}{(T_*-t)^{\frac{1}{\a+1}}}\right),\n \\
 b(x,t)&=&0\n \\
 p(x,t)&=& -\frac{1}{(\a+1)(T_*-t)^{\frac{\a+2}{\a+1}}} (x\cdot
 \nabla ) h\left(\frac{x}{(T_*-t)^{\frac{1}{\a+1}}}\right)\\
 &&\qquad-\frac{1}{2(T_* -t)^{\frac{2}{\a+1}}}
 \left|\nabla
 h\left(\frac{x}{(T_*-t)^{\frac{1}{\a+1}}}\right)\right|^2
 \eq
with a non-constant harmonic function $h$ and $\a\neq -1$, for which
$V=\nabla h \notin C_0 ^1 (\Bbb R^n)$.  Another reason for the
condition (i) is that, since the local classical solution $v$ given
by (\ref{vel}) satisfy $v(\cdot, t)\in C_0 ^1(\Bbb R^n)$ for $t\in
(0, T_*)$, it guarantees the existence of back-to-label map
generated by the velocity field $v$(see \cite{con1}), which is
importantly used in the proof(see
the proof of Theorem 1.2 below).\\
\ \\
\noindent{\textsf{Remark} 1.2} In order to illustrate the decay
condition (ii) we make the following observations. If $\O \in
L^\infty (\Bbb R^n)$ and there exist constants $R, C, \vare_1,
\vare_2
>0$ such that $ |\O (x)|\leq C e^{-\vare_1 |x|^{\vare_2}}$ for $|x|>R$,
then we have
 $\O\in L^q (\Bbb R^n)$ for all $q\in (0, 1)$.
 Indeed,
for all $q\in (0, 1)$, we have
 \bqn
\int_{\Bbb R^n} |\O(x)|^q dx
&=& \int_{|x|\leq R} |\O(x)|^q dx +\int_{|x|>R}  |\O (x)|^q \,dx\\
&\leq&|B_R |^{1-q}\left(\int_{|x|\leq R} |\O(x)| dx\right)^{q} + C^q
\int_{\Bbb R^n} e^{-q\vare_1|x|^{\vare_2}}dx <\infty ,
 \eqn
where $|B_R|$ is the volume of $B_R$.\\
\ \\

Theorem 1.1 will  follow as a corollary of the following more
general theorem.
\begin{theorem}
 Suppose there exists $T>0$ such that we have  a
 representation of the solution $(b, v )$
  to (MHD) by
 \bb\label{thm12a}
 \omega(x,t)=\Psi_1(t) \Omega (\Phi_1 (t)x),
 \ee
 \bb\label{thm12b}
 b(x,t) = \Psi_2(t)B (\Phi_2(t)x)
 \ee
for all $ t\in [0, T)$, where  $\omega= \mathrm{curl}\, v$, $ \Omega
=\mathrm{curl }\, V$ with $V$ satisfying div $V=0$, div $B=0$,
$\Psi_j(\cdot)\in C([0, T );(0, \infty))$,
 $\Phi_j(\cdot)\in C([0, T );\Bbb R^{n\times n})$
 with $\mathrm{det}(\Phi _j (t))\neq 0$ on $[0, T)$ for each $j=1,2$, and $(V, B)$
 satisfying the conditions (i) and (ii) of Theorem 1.1.  Then,
necessarily
 either $\mathrm{det}(\Phi_j(t))\equiv \mathrm{det}(\Phi_j(0))$ on $[0, T)$ for $j=1,2$, or
 $V=B=0$.
\end{theorem}

\subsection{Proof of the theorems}

{\bf Proof of Theorem 1.2} We assume local classical solution $(v,
b)$ of the form (\ref{thm12a})-(\ref{thm12b}). We will show that
this assumption leads to $\Omega=B=0$. By consistency with the
initial condition, $ b_0 (x)=\Psi_2 (0)B(\Phi _2(0)x)$, and hence
  $B (x)=\Psi_2(0)^{-1} b_0 ([\Phi_2 (0) ]^{-1} x)$ for all
$x\in \Bbb R^n$. We can rewrite the representation (\ref{thm12b}) in
the form,
 \bb \label{thm12aa}
 b (x,t) = G(t)b_0 (F(t)x)
 \qquad \forall t\in [0, T),
 \ee
where $G(t)=\Psi_2 (t)/\Psi_2 (0)$, $F(t)=[\Phi _2(0)]^{-1}\Phi_2
(t)$. In order to prove the theorem it suffices to show that either
det$(F(t))=1$ for all $t\in [0,T)$, or $b_0 =0$, since det$(F(t))$=
det$(\Phi _2(t))/$det$(\Phi_2(0))$.
   Let $a\mapsto X(a,t)$ be the
particle
 trajectory mapping, defined by the ordinary differential equations,
 $$
 \frac{\partial X(a,t)}{\partial t} =v(X(a,t),t) \quad;\quad X(a,0)=a.
 $$
We set $A(x,t):=X^{-1} (x ,t)$, which is called the back-to-label
map, satisfying
 \bb\label{inverse}
 A (X(a,t),t)=a, \quad X(A
 (x,t),t)=x.
 \ee
Taking dot product the second equation of (MHD) by $b$, we obtain
 \bb\label{new0}
 \frac{\partial |b|}{\partial t} +(v\cdot \nabla )|b |=\a |b|,
 \ee
where $\a (x,t)$ is defined as
 $$
 \a (x,t)=\left\{ \aligned &\sum_{i,j=1}^n S_{ij}
 (x,t)\xi_i(x,t)\xi_j (x,t) &\quad \mbox{if} \quad b (x,t)\neq 0\\
 &\qquad 0 &\quad \mbox{if} \quad b (x,t)=0\endaligned \right.
 $$
 with
 $$ S_{ij} =\frac12 \left( \frac{\partial v_j}{\partial x_i}
 +\frac{\partial v_i}{\partial x_j}\right)\quad\mbox{and}\quad
\xi (x,t)=\frac{b(x,t)}{|b (x,t)|}.
$$
 In terms of the particle trajectory mapping we
can rewrite (\ref{new0}) as
 \bb\label{new1}
 \frac{\partial }{\partial t} |b
(X(a,t),t)|=\alpha (X(a,t),t) |b (X(a,t),t)|.
 \ee
Integrating (\ref{new1}) along the particle trajectories $\{
X(a,t)\}$, we have
 \bb\label{new2}
 |b (X(a,t),t)|=|b_0 (a)|\exp \left[ \int_0 ^t \a
(X(a,s),s) ds \right].
  \ee
Taking into account the simple estimates
$$ -\|\nabla v(\cdot ,t)\|_{L^\infty}\leq \a
(x,t) \leq \|\nabla v (\cdot ,t)\|_{L^\infty} \quad \forall x\in
\Bbb R^n,
 $$
we obtain from (\ref{new2}) that
 \bqn
 \lefteqn{|b_0 (a)|\exp \left[- \int_0 ^t \|\nabla v
      (\cdot,s)\|_{L^\infty} ds \right]\leq |b (X(a,t),t)|}\hspace{1.in}\\
      && \qquad \leq |b_0 (a)|\exp \left[ \int_0 ^t \|\nabla v
      (\cdot,s)\|_{L^\infty} ds \right],
 \eqn
 which, in terms of the back-to-label map, can be rewritten as
  \bq\label{new3}
 \lefteqn{|b_0 (A(x,t))|\exp \left[- \int_0 ^t \|\nabla v
      (\cdot,s)\|_{L^\infty} ds \right]\leq |b
      (x,t)|}\hspace{1.in}\n\\
     && \leq |b_0 (A(x,t))|\exp \left[ \int_0 ^t \|\nabla v
      (\cdot,s)\|_{L^\infty} ds \right].
 \eq
 Combining this with the self-similar representation formula in (\ref{thm12aa}), we have
 \bq
  \label{new4}
 \lefteqn{|b_0 (A(x,t))|\exp \left[- \int_0 ^t \|\nabla v
      (\cdot,s)\|_{L^\infty} ds \right]\leq G(t) |b _0 (F(t)x)|}\hspace{1.3in}\n \\
      &&\leq
      |b_0 (A(x,t))|\exp \left[ \int_0 ^t \|\nabla v
      (\cdot,s)\|_{L^\infty} ds \right].
 \eq
 Given $q\in (0, q_1)$, computing $L^q(\Bbb R^n)$ norm of the each side of (\ref{new4}),
 we derive
 \bq\label{new4a}
  \lefteqn{\|b_0 \|_{L^q} \exp \left[- \int_0 ^t \|\nabla v
      (\cdot,s)\|_{L^\infty} ds \right]\leq G(t)
      [\mathrm{det}(F(t))]^{-\frac1q}
      \|b _0 \|_{L^q} }\hspace{1.5in}\n \\
      &&\leq
      \|b_0 \|_{L^q}\exp \left[ \int_0 ^t \|\nabla v
      (\cdot,s)\|_{L^\infty} ds \right],
\eq
 where we used the fact $\mathrm{det} (\nabla A (x,t))\equiv 1$.
 Now, suppose $B\neq 0$, which is equivalent to assuming that $b_0 \neq
 0$, then we divide (\ref{new4a}) by $\|b_0 \|_{L^q}$ to obtain
 \bq \label{new5}
 \lefteqn{\exp \left[- \int_0 ^t \|\nabla v
      (\cdot,s)\|_{L^\infty} ds \right]\leq G(t)
      [\mathrm{det}(F(t))]^{-\frac1q}
      }\hspace{1.5in}\n \\
   &&\leq \exp \left[ \int_0 ^t \|\nabla v
      (\cdot,s)\|_{L^\infty} ds \right].
 \eq
 If there exists $t_1\in (0, T)$ such that
 $\mathrm{det}(F(t_1))\neq1$, then either $\mathrm{det}(F(t_1))>1$
 or $\mathrm{det}(F(t_1))<1$. In either case, setting $t=t_1$  and
  passing $q\searrow 0$ in
 (\ref{new5}), we deduce that
 $$\int_0 ^{t_1}\|\nabla v
      (\cdot,s)\|_{L^\infty} ds =\infty.
      $$
      This contradicts with the assumption that the flow is smooth on $(0, T)$, i.e $v\in
      C^1 ([0,T); C^1 _0(\Bbb R^n))$, which is implied by the
by the explicit representation formula (\ref{vel})-(\ref{mag}),
combined with the assumption (i). Hence we need to have
      $B=0$. Setting  $B=0$ in the system (MHD), it reduces to the
     incompressible Euler system, and  the vorticity $\o$  satisfies
 \bb\label{vorticityform}
  \frac{\partial \o }{\partial t} +(v\cdot \nabla )\o =(\o\cdot \nabla )v.
 \ee
 Multiplying (\ref{vorticityform}) by $\o$, we obtain that
$$\frac{\partial |\o | }{\partial t} +(v\cdot \nabla )|\o| =\a |\o|
$$
with the same $\a$ as in (\ref{new0}).  From now on  can repeat the
above argument (\ref{new1})-(\ref{new5}) word by word, and conclude
that $\Omega=$ curl $V=0$ on $\Bbb R^n$. Hence, $V=\nabla h$ for
some scalar function $h\in C^2(\Bbb R^n)$. Since we also have $0=$
div $V$= $\Delta h$, we find that $V=\nabla h$ is harmonic in $\Bbb
R^n$. Combining this
with the hypothesis, $V\in C_0^1 (\Bbb R^n)$,  we obtain that $V=0$. $\square$\\
\ \\
\noindent{\bf Proof of Theorem 1.1} We apply Theorem 1.2 with
$$\Phi_j (t)=(T_* -t)^{-\frac{1}{\a +1}}I,\quad\mbox{and}\quad
 \Psi_j (t) =(T_* -t) ^{-1}
 $$
  for  $j=1,2$, where $I$ is the unit matrix in $\Bbb R^{n\times n}$.
If $\a\neq -1$ and $t\neq 0$, then
 $$ \mathrm{det}(\Phi _j(t))=(T_* -t)^{-\frac{n}{\a +1}}
 \neq T_* ^{-\frac{n}{\a +1}} =\mathrm{det}(\Phi_j (0))
 $$
 for each $j=1,2$.
 Hence, we conclude  that $B=V=0$
by  Theorem 1.2. $\square$\\

\section{Asymptotically self-similar singularities}
 \setcounter{equation}{0}

 In this section we consider scenario of
 more refined possibility of `asymptotically self-similar singularity', which
 means that the local in time smooth solution evolves into a
self-similar profile as the possible singularity time is approached.
The meaning of it will be more clear in the statements of Theorem
2.2. A similar notion was
 considered previously by Giga and Kohn in the context of the
 nonlinear scalar heat equation in \cite{gig}.
 Recently, the author of this paper
 (\cite{cha2}) considered similar notion in the context of 3D Euler and the 3D
 Navier-Stokes equations(see also
 \cite{hou}), and was excluded under suitable conditions on the profile.
We apply similar methods used in \cite{cha2} for our case of (MHD).
In the proof in \cite{cha2} we crucially use a new type of
continuation principle of the local solutions of 3D Euler equations,
where the use of critical homogeneous Besov space
$\dot{B}^{0}_{\infty, 1} (\Bbb R^n)$(see subsection of preliminaries
below for the definition) was essential. For our proof
 we need to establish another continuation principle for our local
 solutions of (MHD), which is Theorem 2.1. For such continuation
 principle we use different critical  Besov space, $\dot{B}^{\frac n2}_{2, 1} (\Bbb
 R^n)$, which is technically important in order to handle the
more complicated nonlinear structure in (MHD) than the case of the
Euler system.

\subsection{Preliminaries}

We follow \cite{tri}(see also \cite{che}). Let $\mathcal{S}$ be the
Schwartz class of rapidly decreasing
 functions. Given $f\in \mathcal{S}$,
 its Fourier transform $\mathcal{F}(f)= \hat{f}$ is defined by
 $$
 \hat{f} (\xi)=\frac{1}{(2\pi )^{n/2}}\int e^{-ix\cdot \xi }
 f(x)dx.
  $$
We consider  $\varphi \in \mathcal{S}$ satisfying
 $\textrm{Supp}\, \hat{\varphi} \subset
 \{\xi \in {\mathbb R}^n \, |\,  \frac12 \leq |\xi|\leq
 2\}$,
 and $\hat{\varphi} (\xi)>0 $ if $\frac12 <|\xi|<2$.
 Setting $\hat{\varphi_j } =\hat{\varphi } (2^{-j} \xi )$ (In other words,
 $\varphi_j (x)=2^{jn} \varphi (2^j x )$.), we can adjust the
 normalization constant in front of $\hat{\varphi}$ so that
  $$
 \sum_{j\in \mathbb{Z}}  \hat{\varphi}_j (\xi )=1\quad \forall \xi \in
 {\mathbb R^n}\setminus \{ 0\}.
 $$
 Given $k\in \mathbb{Z}$,
  we define the function $S_k \in \mathcal{S}$ by its Fourier
 transform
 $$
 \hat{S}_k (\xi )= 1-\sum_{j\geq k+1} \hat{\varphi}_j (\xi ).
 $$
 We observe
 \[
 \mbox{Supp $\hat{\varphi}_j \cap$ Supp $\hat{\varphi}_{j'}=\emptyset$ if
 $|j-j'|\geq 2$}.
 \]
 Let $s\in \mathbb R$, $p,q
 \in [0, \infty]$. Given $f\in \mathcal{S'}$, we denote
 $\Delta_j f=\varphi_j* f$.
 Then the homogeneous Besov semi-norm $\|
f\|_{\dot{B}^{s}_{p,q}}$ is defined by
\[
\| f\|_{\dot{B}^s_{p,q}}=\left\{ \begin{array}{ll} \left[
\sum_{-\infty}^{\infty} 2^{jqs} \| \Delta_j f \|_{L^p}^q
\right]^{\frac{1}{q}}\mbox{ if } q\in [1, \infty) \\
\sup_{j} \left[2^{js} \| \Delta_j f \|_{L^p} \right]\mbox{ if
}q=\infty.
\end{array}\right.
\]
The homogeneous Besov space $\dot{B}^{s}_{p,q}$ is a quasi-normed
space with the quasi-norm given by $\| \cdot \|_{\dot{B}^s_{p,q}}$.
The norm $\|\cdot
 \|_{\dot{B}^s_{p,q}}$  is actually defined up to addition of
 polynomials(namely,
  if $f_1-f_2 $ is a polynomial, then both of $f_1$ and $f_2$ give the same norm),
 and the space  $\dot{B}^s_{p,q} (\Bbb R^n)$ is defined as the
 quotient space of  a class of functions with finite norm, $\|\cdot\|_{\dot{B}^s_{p,q}}$,
 divided by the space of polynomials in $\Bbb R^n$.
For $s>0$  the inhomogeneous Besov space norm $\| f\|_{{B}^s_{p,q}}$
of $f\in \mathcal{S}'$ is defined as $\|f\|_{{B}^s_{p,q}}=\|
f\|_{L^p}+\| f\|_{\dot{B}^s_{p,q}}$. Let us now state some basic
properties for the Besov spaces.
\begin{description}
\item[(i)] Bernstein's Lemma: Assume that $f\in L^p$, $ 1\leq p
\leq \infty$, and $\mbox{supp } \hat{f} \subset \{ 2^{j-2} \leq
|\xi| < 2^{j} \}$, then there exists a constant $C_k$ such that the
following inequality holds:
\[
C_k^{-1} 2^{jk} \| f\|_{L^p} \leq \| D^{k} f\|_{L^p} \leq C_k 2^{jk}
\| f\|_{L^p}.
\]
\item[(ii)] We have the equivalence of norms
\[
\| D^k f \|_{\dot{B}^s_{p,q}} \sim \| f\|_{\dot{B}^{s+k}_{p,q}}.
\]
\item[(iii)] If $s$ satisfies $s \in (-\frac{n}{p}-1, \frac{n}{p}]$, then we
have
 \bb\label{com} \|[u, \Delta_j ] w\|_{L^p} \leq c_j 2^{-j(s+1)} \|
u\|_{\dot{B}^{\frac{n}{p}+1}_{p,1}} \| w\|_{\dot{B}^{s}_{p,1}}
 \ee
with $\sum_{j\in \mathbb{Z}} c_j \leq 1$. In the above, we denote
\[
 [u, \Delta_j ] w=u\Delta_j w-\Delta_j ( uw).
\]
\item[(iv)] We have the following embedding relations for $1\leq p<\infty$:
\bb\label{emb} \dot{B}^{\frac np}_{p,1}(\Bbb R^n)\hookrightarrow C_0
(\Bbb R^n). \ee
\end{description}
In the above (i) and (ii) are standard(see e.g. \cite{tri, che} for
the proofs). The commutator estimate (iii) is proved in \cite{dan},
and (iv) is proved in \cite{bou}.

\subsection{Statement of the  theorems}

In order to prove our main theorem on the asymptotically
self-similar singularities we first establish the following
continuation principle for local classical solution of (MHD), which
is interesting in itself.
\begin{theorem}
Let $p\in [1, \infty)$ and $(v, b) \in [C([0, T);H^m (\Bbb
R^n))]^2$, $m>n/2+1$,  be a classical solution to (MHD). There
exists an absolute constant $\eta
>0$ such that if
 \bb\label{th11}
 \inf_{0\leq t<T} (T-t)
 \left\{\|\o (t)\|_{\dot{B}^{\frac{n}{2}}_{2, 1} }
 +\|j(t)\|_{\dot{B}^{\frac{n}{2}}_{2, 1} }\right\}<\eta
 ,
 \ee
  where $\o=$curl $v$, $j=$curl $b$, then $(v(x,t), b(x,t))$  can be extended to
  a solution of (MHD) in
  $[0, T+\delta]\times \Bbb R^n$,
 and belongs to $C([0, T+\delta]; H^m (\Bbb R^n))$ for some $\delta
 >0$.
\end{theorem}
 \noindent{\textsf{Remark 2.1}} The proof of the local
existence of solutions to (MHD) for $v_0, b_0 \in H^m (\Bbb R^n)$,
$m>n/2+1$, is standard, adapting for example the proof of local
existence theorem for the Euler equations in \cite{kat}(see also
\cite{maj}). The above theorem says that
 if $T_*$ is the first time of singularity, then we have
 the lower estimate of the blow-up rate; there  exists $t_0\in [0,
 T_*)$ such that
 \bb\label{blow}
 \|\o (t)\|_{\dot{B}^{\frac{n}{2}}_{2, 1} }+
 \|j (t)\|_{\dot{B}^{\frac{n}{2}}_{2, 1} } \geq \frac{C}{T_*-t}\quad
 \forall t\in (t_0 , T_*) .
 \ee
\ \\
 As a consequence of the above theorem we can exclude easily the
possibility of `small' self-similar singularities. The proof is in
the next subsection.
\begin{cor}
Let $(v, b)$ be a classical solution to (MHD). Suppose there exist
$\eta
>0, T >0$ and $ t_0\in (0, T) $ such that  we have representation
$$v (x,t)=\frac{1}{(T-t)^{\frac{\a}{\a+1}}} {\bar{V}}
\left(\frac{x}{(T-t)^{\frac{1}{\a+1}}} \right) \quad \forall
(x,t)\in \Bbb R^n \times (t_0, T),
$$
$$b(x,t)=\frac{1}{(T-t)^{\frac{\a}{\a+1}}} {\bar{B}}
\left(\frac{x}{(T-t)^{\frac{1}{\a+1}}} \right) \quad \forall
(x,t)\in \Bbb R^n \times (t_0, T),
$$
where  $(\bar{V}, \bar{B})\in [H^m (\Bbb R^n)]^2$, $m>n/2+1$, and
$\bar{\O}$=curl $\bar{V}$,  $\bar{J}$=curl $\bar{B}$ satisfy
  $$\|\bar{\O}\|_{\dot{B}^{\frac{n}{2}}_{2 ,1}} +
  \|\bar{J}\|_{\dot{B}^{\frac{n}{2}}_{2 ,1}}<\eta .
  $$
  Then,  $\bar{V}=\bar{B}=0$.
\end{cor}
In particular we note that we have milder decay condition  at
infinity for $\bar{\Omega}$ in the above corollary  than in Theorem
1.1. \\
\ \\
The following theorem exclude the possibility of a type of
asymptotically self-similar singularity for (MHD).
\begin{theorem}
Let $(v, b) \in [C([0, T);H^m (\Bbb R^n))]^2$, $m>n/2+1$,  be a
classical solutions to (MHD).
 Suppose there exist  functions $\bar{V}, \bar{B}$ satisfying the condition (ii) for
 $V, B$ in Theorem 1.1 such that the following boundedness and the convergence hold true:
 \bq\label{th122}
 \lefteqn{\sup_{0<t<T}
(T-t)^{\frac{\a-n}{\a+1}}\left\|v(\cdot, t)
-\frac{1}{(T-t)^{\frac{\a}{\a+1}}} \bar{V}
\left(\frac{\cdot}{(T-t)^{\frac{1}{\a+1}}} \right)\right\|_{L^1
}}\hspace{.0in}\n \\
&&+ \sup_{0<t< T} (T-t)^{\frac{\a-n}{\a+1}}\left\|b(\cdot, t)
-\frac{1}{(T-t)^{\frac{\a}{\a+1}}} \bar{B}
\left(\frac{\cdot}{(T-t)^{\frac{1}{\a+1}}} \right)\right\|_{L^1 }
<\infty,\n \\
  \eq
  and
  \bq\label{th12}
   \lefteqn{\lim_{t\nearrow T}
(T-t) \left\|\o(\cdot, t) -\frac{1}{T-t} \bar{\O}
\left(\frac{\cdot}{(T-t)^{\frac{1}{\a+1}}}
\right)\right\|_{\dot{B}^{\frac{n}{2}}_{2, 1} }}\hspace{.0in}\n \\
 &&+
 \lim_{t\nearrow T}
(T-t)\left\|j(\cdot, t) -\frac{1}{T-t}  \bar{J}
\left(\frac{\cdot}{(T-t)^{\frac{1}{\a+1}}}
\right)\right\|_{\dot{B}^{\frac{n}{2}}_{2, 1} } =0,\n \\
  \eq
  where $\bar{\O}$=curl $\bar{V}$ and $\bar{J}$=curl
$\bar{B}$.
 Then, $\bar{V}=\bar{B}=0$, and $(v, b)$ can be extended to a solution of
  (MHD) in $[0, T+\delta]\times \Bbb R^n$,
 and belongs to $C([0, T+\delta]; H^m (\Bbb R^n))$ for some $\delta
 >0$.
\end{theorem}
\noindent{\textsf{Remark 2.3}} We note that we did not impose the
condition (i) for $\bar{V}$ and $\bar{B}$. This condition is
guaranteed by
the assumption of (\ref{th122}) and (\ref{th12}).\\
\ \\
\noindent{\textsf{Remark 2.4}} We note that Theorem 1.2 still
 does not exclude the possibility that the vorticity and the magnetic field
  converge to the
 asymptotically self-similar singularity in the weaker sense than $L^\infty$
 sense. Namely, a self-similar vorticity profile could be approached
 from a local classical solution in the pointwise sense in space,
  or in the $L^q(\Bbb R^n)$ sense for
some $q\in (1, \infty)$.\\
 \ \\
As an immediate corollary of Theorem 2.2 we have the following
information of the behaviors of solution near possible singularity,
which is not necessarily of the self-similar type.
\begin{cor}
Let $(v, b) \in [C([0, T_*);H^m (\Bbb R^n))]^2$, $m>n/2+1$,  be a
classical solutions to (MHD), which blows up at $T_*$.  We expand
the solution of the form:
 \bq
 v(x,t)&=&\frac{1}{(T_*-t)^{\frac{\a}{\a+1}}}
\bar{V}\left(\frac{x}{(T_*-t)^{\frac{1}{\a+1}}}\right) +\bar{v}(x,t),\\
b(x,t)&=&\frac{1}{(T_*-t)^{\frac{\a}{\a+1}}}
\bar{B}\left(\frac{x}{(T_*-t)^{\frac{1}{\a+1}}}\right)
+\bar{b}(x,t),
 \eq
 where $(\bar{V}, \bar{B})$ satisfies the conditions (i)-(ii) for $(V,B)$ in Theorem
 1.1, and $\a \neq -1$, $t\in [0, T_*)$.
 Then, either
 \bb
 \lim\sup_{t\nearrow T_*} (T_*-t)^{\frac{\a-n}{\a+1}}\|\bar{v}(\cdot, t)\|_{L^1 }=\infty,
 \ee
 or there exists $\e_0 >0$ such that
 \bb
 \lim\sup_{t\nearrow T_*}(T_*-t)
 \left(\|
 \bar{\o}(t)\|_{\dot{B}^{\frac{n}{2}}_{2, 1} }
 +\| \bar{j}(t)\|_{\dot{B}^{\frac{n}{2}}_{2, 1} }\right) >\e_0,
 \ee
 where $\bar{\o}=$curl $\bar{v}$, and $\bar{j}=$curl $\bar{b}$
 respectively.
\end{cor}
\subsection{Proof of the theorems}
\noindent{\bf Proof of Theorem 2.1}
  We set
  $$
  v_\pm=v\pm b, \quad \pi =p+\frac12 |b|^2.
 $$
  Adding and subtracting the first and the second equations of (MHD), we obtain
  \bb\label{mhdpm}
 \frac{\partial v_\pm}{\partial t} +(v_\mp\cdot \nabla )v_\pm=-\nabla
 \pi.
 \ee
 Taking the operation $\Delta_j $ on (\ref{mhdpm}), we have
\begin{equation}\label{ns1}
\frac{\partial}{\partial t} \Delta_j v_\pm +  (v_\mp\cdot \nabla
)\Delta_j v_\pm =[v_\mp, \Delta_j ]\cdot \nabla v_\pm+\nabla
(\Delta_j \pi) .
\end{equation}
Multiplying $\Delta_j v_\pm $ on the both sides of (\ref{ns1}) and
integrating over $\mathbb{R}^n$, we obtain, after integration by
part,
\begin{equation}\label{ns20}
\frac12 \frac{d}{dt} \| \Delta_j v_\pm\|_{L^2}^2  \leq \| [\Delta_j,
v_\mp]\cdot \nabla v_\pm \|_{L^2} \| \Delta_j v_\pm \|_{L^2}.
\end{equation}
Dividing the both sides of (\ref{ns20}) by $\| \Delta_j v_\pm
\|_{L^2}$, and using the commutator estimates (\ref{com}), we deduce
\begin{equation}\label{ns3}
\frac{d}{dt} \| \Delta_j v_\pm \|_{L^2} \leq C c_j
2^{-(\frac{n}{2}+1)j} \left(\| v_\pm
\|_{\dot{B}^{\frac{n}{2}+1}_{2,1}} +\|
v_\mp\|_{\dot{B}^{\frac{n}{2}+1}_{2,1}}\right) \| v_\pm
\|_{\dot{B}^{\frac{n}{2}+1}_{2,1}}.
\end{equation}
Multiplying $2^{(\frac{n}{2}+1)j}$ on the both sides of (\ref{ns3}),
and taking summation over $j\in \mathbb{Z}$, it reduces that
 \bb\label{mhd1}
\frac{d}{dt} \| v_\pm (t) \|_{\dot{B}^{\frac{n}{2}+1}_{2,1}}
 \leq C \left(\|
v_\pm (t)\|_{\dot{B}^{\frac{n}{2}+1}_{2,1}} +\| v_\mp
(t)\|_{\dot{B}^{\frac{n}{2}+1}_{2,1}}\right) \| v_\pm (t)
\|_{\dot{B}^{\frac{n}{2}+1}_{2,1}}.
 \ee
  Let us define
 $$
 X(t)=\|v_+ (t)\|_{\dot{B}^{\frac{n}{2}+1}_{2,1}}
 +\|v_- (t)\|_{\dot{B}^{\frac{n}{2}+1}_{2,1}}.
 $$
 Then, from (\ref{mhd1}) we have
 \bb
\frac{d }{dt}X(t) \leq C X(t)^2.
 \ee
By the Gronwall lemma we deduce
 \bb\label{gro} X(t) \leq
\frac{X(0)}{1-C t X(0)}.
  \ee
  Translating in time, we have instead of (\ref{gro})
 \bb\label{pro11}
X(T_1) \leq \frac{C_1 X(t)}{1-C_2(T_1-t)X(t)}
 \ee
 for all $T_1 \in (0,T)$ and $t\in [0, T_1)$, where   $C_1, C_2$ are absolute
 constants.
We set $\eta=\frac{1}{2C_2}$. For such $\eta$, we  suppose
(\ref{th11}) holds true. Then,  there exists $t_1 \in  [0, T)$ such
that $(T-t_1) X(t_1)<\eta$. Fixing $t=t_1$, and passing $T_1
\nearrow T$ in (\ref{pro11}), we find that
$$\lim\sup_{T_1 \nearrow T} X(T_1 ) \leq  2C_1 X( t_1 )
<\infty ,$$
 and hence,  we obtain the estimate for the homogeneous Besov norm,
 \bb\label{local1}
 \sup_{0< t < T}
 \left(\|v(t)\|_{\dot{B}^{\frac{n}{2}+1}_{2,1}} +
 \| b(t)\|_{\dot{B}^{\frac{n}{2}+1}_{2,1}}
 \right) < C \sup_{0< t< T} X(t):=  M(T) <\infty.
\ee
 if our hypothesis (\ref{th11}) holds true.
Therefore,  we have
 \bqn
  \lefteqn{\int_0 ^T\left\{ \|\o (t)\|_{L^\infty}
  +\|j (t)\|_{L^\infty}\right\}dt\leq \int_0 ^T\left\{ \|\nabla v(t)\|_{L^\infty}
  +\|\nabla b(t)\|_{L^\infty}\right\}dt}\hspace{.1in}\\
 &&\leq C\int_0 ^T \left\{\| v(t)\|_{\dot{B}^{\frac{n}{2}+1}_{2,1}}
  +\|b(t)\|_{\dot{B}^{\frac{n}{2}+1}_{2,1}}\right\}dt  \leq C M(T) T <\infty,
  \eqn
 where we used the embedding $\dot{B}^{\frac{n}{2}}_{2,1} (\Bbb R^n) \hookrightarrow
 L^\infty (\Bbb R^n )$.
Applying the BKM type of blow-up criterion for (MHD) derived in
\cite{caf}, we can continue our local solution
  $(v(x, t), b(x, t))$ until $t=T+\delta $, and
  $(v, b)\in [C([0, T+\delta];H^m (\Bbb R^n ))]^2$ for some $\delta >0$, where $m>n/2+1$. $\square$\\
\ \\
\noindent{\bf Proof of Corollary 2.1} We just observe that
$$
 (T-t)\|\o (t)\|_{\dot{B}^{\frac{n}{2}}_{2,1}}=
\|\bar{\O}\|_{\dot{B}^{\frac{n}{2}}_{2,1}},\quad
  (T-t)\|j (t)\|_{\dot{B}^{\frac{n}{2}}_{2,1}}=
\|\bar{J}\|_{\dot{B}^{\frac{n}{2}}_{2,1}}
 $$
for all $t\in (t_0, T)$. Hence, our smallness condition, $
 \|\bar{\Omega}\|_{\dot{B}^{\frac{n}{2}}_{2,1}}
 +\|\bar{J}\|_{\dot{B}^{\frac{n}{2}}_{2,1}}<\eta,
 $
leads to
$$\inf_{t_0<t<T}(T-t)\left\{ \|\o (t)\|_{\dot{B}^{\frac{n}{2}}_{2,1}}
 +\|j (t)\|_{\dot{B}^{\frac{n}{2}}_{2,1}}\right\}
<\eta.
 $$
 Applying Theorem 2.1, for initial time  at $t=t_0$, we
conclude  that $(v, b)\in [C^1 ([t_0, T); H^m (\Bbb R^n ))]^2$
cannot have singularity at $t=t_0$, hence we need to have $V=B=0$. $\square$\\
\ \\
\noindent{\bf Proof of Theorem 2.2} We change variables from the
physical ones $(x,t) \in \Bbb R^n \times [0,T)$ to the `self-similar
variables' $(y,s)\in \Bbb R^n\times [0, \infty)$ as follows:
$$
y=\frac{x}{(T-t)^{\frac{1}{\a+1}}}, \quad s=\frac{1}{\a+1} \log
\left( \frac{T}{T-t}\right).
$$
Based on this change of variables, we transform the functions
$(v,p)\mapsto (V, P)$ according to
 \bq
 v(x,t)&=&\frac{1}{(T-t)^\frac{\a}{\a+1}} V(y,s ),\\
 b(x,t)&=&\frac{1}{(T-t)^\frac{\a}{\a+1}} B(y,s ),\\
 p (x,t)&=&\frac{1}{(T-t)^\frac{2\a}{\a+1}}P (y,s ).
 \eq
 Substituting $(v,b, p )$ into the $(MHD)$, we obtain the
 following equivalent evolution equations for
 $(V,P)$,
 \bb\label{selfmhd}
 \left\{ \aligned
  &  \frac{1 }{\a+1}V_s +\frac{\a}{\a+1} V +\frac{1}{\a +1}(y \cdot \nabla)V
  + (V\cdot \nabla )V =(B\cdot \nabla )B-\nabla
( P +\frac12 |V|^2),\\
  &\frac{1 }{\a+1}B_s +\frac{\a}{\a+1} B +\frac{1}{\a +1}(y \cdot
  \nabla)B
  + (V\cdot \nabla )B =(B\cdot \nabla )V,\\
 & \qquad\mathrm{div}\, V=\mathrm{div}\, B=0,\\
 & V(y,0)=V_0 (y)=T^{\frac{\a}{\a +1}} v_0 (T^{\frac{1}{\a +1}}y),\quad
  B(y,0)=B_0 (y)=T^{\frac{\a}{\a +1}} b_0 (T^{\frac{1}{\a +1}}y).
  \endaligned \right.
 \ee
 In terms of $(V,B)$ the conditions
 (\ref{th122}) and
 (\ref{th12}) are translated
into
 \bb\label{econv11}
 \sup_{0<s <\infty}\|V(\cdot ,s) -\bar{V}(\cdot )\|_{L^1}+ \sup_{0<s< \infty}\|B(\cdot ,s) -\bar{B}(\cdot )\|_{L^1}<\infty,
 \ee
 and
 \bb\label{econv1}
 \lim_{s\to \infty}\|\O(\cdot ,s) -\bar{\O}(\cdot )\|_{\dot{B}^{\frac{n}{2}}_{2, 1}}=
 \lim_{s\to \infty}\|J(\cdot ,s) -\bar{J}(\cdot )\|_{\dot{B}^{\frac{n}{2}}_{2,
 1}}=0.
 \ee
From the fact that the Calderon-Zygmund singular integral operator
maps
 $\dot{B}^{\frac{n}{2}}_{2, 1} (\Bbb R^n)$ into
 itself boundedly, we obtain from (\ref{econv1})
  \bb\label{econv1a}
 \lim_{s\to \infty}\|V(\cdot ,s) -\bar{V}(\cdot )\|_{\dot{B}^{\frac{n}{2}+1}_{2, 1}}=
 \lim_{s\to \infty}\|B(\cdot ,s) -\bar{B}(\cdot )\|_{\dot{B}^{\frac{n}{2}+1}_{2,
 1}}=0,
 \ee
 from which, thanks to the embedding (\ref{emb})   and (\ref{econv11}),
 we have
 \bb\label{econv2}
  \lim_{s\to \infty}\| V(\cdot ,s) -\bar{V}(\cdot )\|_{C^1(B_R) }=
   \lim_{s\to \infty}\| B(\cdot ,s) -\bar{B}(\cdot
   )\|_{C^1(B_R)}=0
  \ee
  for all $R>0$, where  $B_R=\{ x\in \Bbb R^n \, |\, |x|<R\}$.
  Moreover, we find that $\bar{V}, \bar{B}\in C_0 (\Bbb R^n)$,  satisfying the condition
  (i) of Theorem 1.1.
 Similarly to \cite{hou}, we consider scalar test functions
$\xi \in C^1_c (0,1)$ with $\int_0 ^1\xi (s)ds\neq 0$, $\psi \in
C^1_c (\Bbb R^n )$ and the vector test function $\phi =(\phi_1 ,
\cdots, \phi_n )\in C_c^1 (\Bbb
 R^n )$ with div $\phi=0$. We multiply the first equation of
 (\ref{selfmhd}) by $\xi
 (s-k)\phi (y)$, and integrate it over $\Bbb R^n\times [k, k+1]$,
 and then we integrate by part for the terms including the
 time derivative and the pressure term to obtain
 \bq\label{keq1}
 &&-\int_0^{1}\int_{\Bbb R^n} \xi _s(s) \phi(y)\cdot V(y,s+k)
 dyds\n\\
 &&+\int_0 ^{1}\int_{\Bbb R^n}\xi (s)\phi(y) \cdot[\a V +(y \cdot \nabla)V +(\a+1)
 (V\cdot \nabla )V](y,s+k)  dyds\n\\
 &&\quad -(\a+1)\int_0 ^{1}\int_{\Bbb R^n}\xi (s)\phi(y) \cdot(B\cdot \nabla )B](y,s+k)
 dyds=0,
 \eq
 and
\bq\label{keq2}
  &&-\int_0^{1}\int_{\Bbb R^n} \xi _s(s) \psi(y) B(y,s+k)
 dyds\n\\
 &&+\int_0 ^{1}\int_{\Bbb R^n}\xi (s)\psi(y) [\a B +(y \cdot \nabla)B +(\a+1)
 (V\cdot \nabla )B](y,s+k)  dyds\n \\
 &&\quad -(\a+1)\int_0 ^{1}\int_{\Bbb R^n}\xi (s)\psi(y)(B\cdot \nabla )V](y,s+k)  dyds=0.
 \eq
Passing to the limit $k\to \infty$ in (\ref{keq1})-(\ref{keq2}),
using the
  convergence (\ref{econv2}), $\int_0 ^1\xi _s(s)ds=0$ and $\int_0 ^1\xi
  (s)ds\neq 0$,
  we  find that $\bar{V}, \bar{B}\in C^1 _0(\Bbb R^n)$ satisfies
$$
\int_{\Bbb R^n} [\a \bar{V} +(y \cdot \nabla)\bar{V} +(\a+1)
(\bar{V}\cdot \nabla
 )\bar{V} -(\a+1)
(\bar{B}\cdot \nabla
 )\bar{B}]\cdot \phi dy=0,
 $$
 and
 $$
\int_{\Bbb R^n} [\a \bar{B} +(y \cdot \nabla)\bar{B} +(\a+1)
(\bar{V}\cdot \nabla
 )\bar{B} -(\a+1)
(\bar{B}\cdot \nabla
 )\bar{V}]\psi dy=0,
 $$
 for all vector test function $\phi \in C_c^1 (\Bbb
 R^n)$ with div $\phi=0$, and scalar test function $\psi \in C^1_c
(\Bbb R^n )$
 Hence, there exists a scalar function
 $\bar{P'}$, which can be written without loss of generality that $\bar{P'}=
 \bar{P}+ \frac12 |\bar{B}|^2$ for another scalar function $\bar{P}$,  such that
 \bq\label{eleray1}
\lefteqn{\a  \bar{V}+ (y\cdot\nabla )\bar{ V}
 +(\a+1)(\bar{V}\cdot \nabla )\bar{V} =(\a+1)(\bar{B}\cdot \nabla )\bar{B}}\hspace{2.in}\n\\
 &&\qquad \quad-(\a+1)\nabla
(\bar{P}+ \frac12 |\bar{B}|^2),
 \eq
 and
 \bb\label{eleray2}
  \a  \bar{B}+ (y\cdot\nabla )\bar{ B}
 +(\a+1)(\bar{V}\cdot \nabla )\bar{B} =(\a+1)(\bar{B}\cdot \nabla
 )\bar{V}.
 \ee
 On the other hand, we can pass $s\to \infty$ directly in the
incompressibility equations for $V$ and $B$ in (\ref{selfmhd}) to
have
 \bb\label{eleray3}
\mathrm{ div} \, \bar{V}=\mathrm{div}\, \bar{B}=0.
 \ee
The equations (\ref{eleray1})-(\ref{eleray3})  show that $(\bar{V},
\bar{B} )$ is a classical solution of (\ref{mhdleray}).  Since, by
hypothesis, curl $\bar{V}=\bar{\Omega}$ and $\bar{B}$ satisfy the
condition (ii) of Theorem 1.1, we can deduce $\bar{V}=\bar{B}=0$ by
that theorem.
 Hence,  (\ref{econv1}) leads to
 $$\lim_{s\to \infty} \|\O
 (s)\|_{\dot{B}^{\frac{n}{2}}_{2, 1}}=\lim_{s\to \infty} \|J
 (s)\|_{\dot{B}^{\frac{n}{2}}_{2, 1}}=0.
 $$
 Thus, for $\eta >0$ given in Theorem 1.1,
 there exists $s_1>0$ such that
 $$\|
 \O (s_1 )\|_{\dot{B}^{\frac{n}{2}}_{2, 1}}
 +\|J  (s_1 )\|_{\dot{B}^{\frac{n}{2}}_{2, 1}}<\eta .
 $$
Let us set $t_1=T[1-e^{(\a +1)s_1} ]$. Going back to the original
physical variables, we have
 $$
  (T-t_1)\|\o (t_1)\|_{\dot{B}^{\frac{n}{2}}_{2, 1}} +
  (T-t_1)\|j (t_1)\|_{\dot{B}^{\frac{n}{2}}_{2, 1}}<\eta .
 $$
 Applying Theorem 2.1, we conclude the proof. $\square$

 \section{The 2D magnetohydrodynamics}
 \setcounter{equation}{0}

 In the 2D magnetohydrodynamics, namely (MHD) with $n=2$, we have the following reduced form of (MHD),
 using the stream function for the magnetic field, $b=\nabla ^\bot \varphi$.
\[
\mathrm{ (2D-MHD)}
 \left\{ \aligned
 &\frac{\partial v}{\partial t} +(v\cdot \nabla )v =(b\cdot\nabla)b-\nabla (p +\frac12 |b|^2),\\
 &\frac{\partial \varphi }{\partial t} +(v\cdot \nabla )\varphi =0, \\
 &\quad \mathrm{div}\, v=0, \quad b=\nabla^\bot\varphi, \\
  &\varphi(x,0)=\varphi_0 (x), \quad v(x,0)=v_0 (x)
  \endaligned
  \right.
  \]
  Unlike the case of 2D Euler equations, the question
  of finite time blow-up/global regularity question  in the magnetohydrodynamics
  is wide open even for the 2D case. Similar situation is for the
  question of global existence of weak solutions.
  We mention the result of
  global existence of weak solutions for a partial viscosity case obtained by Kozono in
  \cite{koz}. There are studies on the
  possible scenarios of finite time singularity other than that of self-similar type in
  \cite{cor1, cor2}. We note that the equation for $\varphi$
  in (2D-MHD) is nothing but a
  divergence free transport equation, for which the nonexistence of
  self-similar singularity is obtained in \cite{cha1} in a general setting.
Applying the general theorem  we first have the following  improved
version of Theorem 1.2.
\begin{theorem}
 Suppose there exists $T>0$ such that we have  a
 representation of the solution $(v, b)$
  to the system (2D-MHD)  with  $b=\nabla^\bot
  \varphi$ such that we have the representation:
 \bb
 \varphi (x,t)=\Psi(t) \Theta(\Phi (t)x),
 \ee
for all $ t\in [0, T)$, where $\Psi(\cdot)\in C([0, T );(0,
\infty))$,
 $\Phi(\cdot)\in C([0, T );\Bbb R^{2\times 2})$
 with $\mathrm{det}(\Phi  (t))\neq 0$ on $[0, T)$. Suppose also that there exist
 $q_1<q_2$ with  $q_1, q_2 \in (0, \infty]$ such that $\Theta \in L^{q_1}
 (\Bbb R^2 )\cap L^{q_2} (\Bbb R^2 )$,  then either det$(\Phi (t))=$det$(\Phi (0))$
 or  $\Theta =0$.
\end{theorem}
The proof follows immediately by Theorem 2.2 in \cite{cha1} to the
first equation of (2D-MHD). Note that we do not assume any
integrability condition for the vorticity in the above theorem. As a
corollary of the above theorem we obtain:
\begin{theorem}
Suppose there exists $T_*>0$ such that we have  a
 solution $(v, b)$
 with $b=\nabla^\bot \varphi$ to the system (2D-MHD) with the representation:
  \bq\label{th31}
  v(x,t)&=&\frac{1}{(T_*-t)^{\frac{\a }{\a+1}}}
V \left(\frac{x}{(T_*-t)^{\frac{1}{\a+1}}}\right)\\
\label{th32}
  \varphi (x, t)&=&\frac{1}{(T_*-t)^{\frac{\a -1}{\a+1}}}
\Theta\left(\frac{x}{(T_*-t)^{\frac{1}{\a+1}}}\right)
 \eq
 for $t \in (0 ,T_*)$ with $\a \neq -1$, and let $V\in H^m (\Bbb R^2)$,
 $m>2$, and  $\Theta \in L^{q_1}
 (\Bbb R^2 )\cap L^{q_2} (\Bbb R^2 )$
 for some  $q_1, q_2 \in (0, \infty]$  with $q_1<q_2$. Then, $\Theta=0$,
  and  $V=0$.
\end{theorem}
\noindent{\textsf{Remark} 3.1} We note that the integrability
condition on the vorticity $\O=$ curl $V\in H^{m-1} (\Bbb R^2)$,
$m>2$, and the corresponding decay at infinity
 are weaker here than in Theorem 1.1. The assumption for
$\Theta$ in the above theorem is satisfied if
 $$
B=\nabla^\bot\Theta \in W^{-1,q_1}
 (\Bbb R^2 )\cap W^{-1,q_2} (\Bbb R^2 )\quad \forall \,1< q_1< q_2 \leq \infty.
  $$
  Since $L^{\frac{2q}{2-q}}(\Bbb R^2)\hookrightarrow W^{-1,q} (\Bbb
  R^2 )$ for $q\in (1, 2)$, we have that the representation (\ref{vel})-(\ref{mag})  for the solution $(v,b)$
  of (MHD) for $n=2$ with
  $B$  satisfying
  $$B\in L^{p_1} (\Bbb R^2 )\cap L^{p_2} (\Bbb R^2), \quad \forall\,
  2< p_1<p_2 <\infty
  $$
  implies $V=B=0$.\\
  \ \\
  \noindent{\bf Proof of Theorem 3.2} Applying Theorem 3.1 to the
  representation of $\varphi$ in (\ref{th32}), we first conclude that
  $\Theta=\varphi=0$. Setting  $\varphi =0$ in (2D-MHD), then the system reduces to
  the 2D Euler equations, for which we have the  global well-posedness of the classical
  solutions for initial data  $v(t) \in H^m(\Bbb R^2)$, $m>2$, at $t\in [0, T_*)$, and
  hence the representation for $v$ in (\ref{th31}) is possible only when
  $V=0$. $\square$\\
  \ \\
Combining Theorem 3.2 with Theorem 2.2, we also immediately obtain
the following theorem.
\begin{theorem}
Let $(v, b) \in [C([0, T);H^m (\Bbb R^2))]^2$, $m>2$,  be a
classical solutions to (MHD) with $n=2$.
 Suppose there exist  functions $\bar{V}, \bar{B}=\nabla^\bot \bar{\Theta}$
 with $\bar{V}$
 and $\bar{\Theta}$,
  satisfying the conditions of $V$ and $\Theta$ respectively in Theorem
  3.2 and $\a \neq -1$
   such that the following estimate and convergence hold true:
 \bq\label{th122a}
 \lefteqn{\sup_{0<t< T}
(T-t)^{\frac{\a-2}{\a+1}}\left\|v(\cdot, t)
-\frac{1}{(T-t)^{\frac{\a}{\a+1}}} \bar{V}
\left(\frac{\cdot}{(T-t)^{\frac{1}{\a+1}}} \right)\right\|_{L^1
}}\hspace{.0in}\n \\
&&+ \sup_{0<t< T}(T-t)^{\frac{\a-2}{\a+1}}\left\|b(\cdot, t)
-\frac{1}{(T-t)^{\frac{\a}{\a+1}}} \bar{B}
\left(\frac{\cdot}{(T-t)^{\frac{1}{\a+1}}} \right)\right\|_{L^1 }
<\infty,\n \\
  \eq
  and
  \bq\label{th12b}
   \lefteqn{\lim_{t\nearrow T}
(T-t) \left\|\o(\cdot, t) -\frac{1}{T-t} \bar{\O}
\left(\frac{\cdot}{(T-t)^{\frac{1}{\a+1}}}
\right)\right\|_{\dot{B}^{1}_{2, 1} }}\hspace{.0in}\n \\
 &&+ \lim_{t\nearrow T}
(T-t)\left\|j(\cdot, t) -\frac{1}{T-t}  \bar{J}
\left(\frac{\cdot}{(T-t)^{\frac{1}{\a+1}}}
\right)\right\|_{\dot{B}^{1}_{2, 1} } =0,\n \\
  \eq
  where we set
  $$\o=\partial_{1} v_2-\partial_{2} v_1, \quad j=\partial_{1} b_2-\partial_{2}
  b_1 ,
  $$
  and
  $$
  \bar{\O}=\partial_{1}\bar{ V}_2-\partial_{2} \bar{V}_1, \quad \bar{J}=\Delta
 \bar{ \Theta}.
  $$
 Then, $\bar{V}=\bar{B}=0$, and $(v, b)$ can be extended to a solution of
 (MHD) in $[0, T+\delta]\times \Bbb R^2$,
 and belongs to $C([0, T+\delta]; H^m (\Bbb R^2))$, $m>2$,  for some $\delta
 >0$.
\end{theorem}

\end{document}